\documentclass[sn-mathphys]{sn-jnl}

\usepackage{amssymb}  

\usepackage{amsthm}
\usepackage{commath}

\usepackage{subcaption}
\usepackage{framed}
\usepackage[normalem]{ulem}
\usepackage{enumerate}
\usepackage{cleveref}
\usepackage{rotating}

\theoremstyle{definition}

\newcommand{\of}{\textsc{OpenFOAM}\textsuperscript\textregistered}

\newcommand{\vel}{\mathbf{v}}
\newcommand{\xb}{\mathbf{x}}
\newcommand{\gb}{\mathbf{g}}
\newcommand{\por}{\phi}

\newcommand{\perm}{\mathbb{K}}

\newcommand{\prgh}{p_{\rho gh}}

\newcommand{\be}[1]{
\begin{equation}
\expandafter\label{eq:#1}
}
\newcommand{\ben}{
\begin{equation}
\expandafter
}
\newcommand{\ee}{\end{equation}}

\newcommand{\bg}[1]{
\begin{gather}\nonumber
\expandafter\label{eq:#1}
}
\newcommand{\eg}{\end{gather}}

\newcommand{\bfig}[1]{
\begin{figure}
\expandafter\label{fig:#1}
}
\newcommand{\efig}{\end{equation}}

\newcommand{\diverg}{\nabla\cdot}

\newcommand{\grad}{\nabla}

\jyear{2022}%

\title
[Computational framework for flow and transport in porous media]
{
Computational framework for complex flow and transport in heterogeneous porous media
}

\author*[1]{\fnm{Matteo} \sur{Icardi}}\email{matteo.icardi@nottingham.ac.uk}

\author[1]{\fnm{Eugenio} \sur{Pescimoro}}

\author[2]{\fnm{Federico} \sur{Municchi}}

\author[3]{\fnm{Juan H.} \sur{Hidalgo}}

\affil[1]{\orgdiv{School of Mathematical Sciences}, \orgname{University of Nottingham}, \orgaddress{\street{University Park}, \city{Nottingham}, \postcode{NG7 2RD}, \state{England}, \country{UK}}}

\affil[2]{\orgname{Colorado School of Mines}, \orgaddress{\street{1500 Illinois St.}, \city{Golden}, \postcode{CO 80401}, \state{Colorado}, \country{USA}}}

\affil[3]{\orgname{Institute of Environmental Assessment and Water Research}, \orgaddress{\street{C. Jordi Girona 18-26}, \city{Barcelona}, \postcode{08034}, \country{Spain}}}

\begin{document}



\abstract{
We present a flexible scalable open-source computational framework, named \texttt{SECUReFoam}, based on the finite-volume library \of, for flow and transport problems in highly heterogeneous geological media and other porous materials. The framework combines geostatistical pre- and post-processing tools with specialised Partial Differential Equations solvers. Random fields, for permeability and other physical properties, are generated by means of continuous or thresholded Gaussian random fields with various covariance/variogram functions. The generation process is based on an explicit spectral Fourier decomposition of the field which, although more computationally intensive than Fast Fourier Transform methods, allows a more flexible choice of statistical parameters and can be used for general geometries and grids. Flow and transport equations are solved for single-phase and variable density problems, with and without the Boussinesq approximation, and for a wide range of density, viscosity, and dispersion models, including dual-continuum (dual permeability or dual porosity) formulations. The mathematical models are here presented in details and the numerical strategies to deal with heterogeneities, equation coupling, and boundary conditions are discussed and benchmarked for the heterogeneous Henry and Horton-Rodgers-Lapwood problems, and other test cases. We show that our framework is capable of dealing with large permeability variances, viscous instabilities, and large-scale three-dimensional transport problems.
}


	
	\maketitle

	\section{Introduction}
        Partial differential equations (PDEs) models in porous media are an essential part of many engineering, environmental and biological applications. They play a fundamental role in geothermal energy utilization \cite{Fan2007, Hidalgo2009,Limberger2018}, aquifer remediation \cite{Zhang2017}, drug delivery in tissues \cite{Khanafer2006} or bio-film formation \cite{Rittmann1993,Gaebler2018}, to name a few. Fluid flow, energy and solute transport in porous media are complex because of the variation in fluid properties and the heterogeneous nature of the porous medium. Fluid properties such as density and viscosity vary with the solute concentration and fluid temperature. These changes can lead to fluid instabilities greatly affecting the migration of dissolved substances, the mixing of fluids and chemical reactions \cite{Homsy1987,Kueper1991, Szulczewski2013, Hidalgo2015}. Heterogeneity is present across all scales. The interaction between the porous structure and the fluid flow manifest itself in the formation of preferential flow paths and stagnant regions. This in turn affects the migration and residence times of solutes, which often display an anomalous behavior \cite{Berkowitz2006}. 

        Given the complexity of the phenomena related to porous media, numerical simulations provide the possibility to explore a multitude of configurations that can help analysing field and experimental data. However, producing accurate and reliable numerical data requires, considering the current state-of-the-art, advanced knowledge of numerical schemes and programming languages. Therefore, there is a need for robust and precise simulation tools capable of providing high quality results, but also flexible enough to deal with the enormous variety of porous media applications. Reproducibility of numerical results is often as important as accuracy. Reproducibility can only be achieved if data sets, results and simulation tools are accessible to the public. This can be achieved through the use of an open-source license and by building the simulation tools over well-known, well-maintained open-source libraries.

In this work we present the open-source computational framework \texttt{SECUReFoam} for the simulation of complex flow and transport in porous media. The framework is developed in \texttt{C++} within the \of\ library and includes new solvers, boundary conditions, pre- and post-processing utilities, and new numerical schemes. We detail the implementation of the geostatistical pre-processing library and solution strategites for saturated flow and transport solvers that account for variable fluid properties (e.g., density, viscosity, relative permeability). Heterogeneity can be accounted directly through the definition of random permeability fields and by means of double porosity models. The formers are based on multi-Gaussian and thresholded Gaussian random fields (also known as pluri-Gaussian) \cite{beucher2016truncated, journel1978mining}. Although this kind of fields can be generated with other available geostatistics toolboxes (GSLib \cite{deutsch1998gslib}, T-PROGS \cite{carle1999t}), the integration with the flow and transport solvers provides the computational framework the capability to tackle a wide variety  complex problems.

Various other open-source packages are available for porous media applications, such as OpenGeoSys \cite{kolditz2012opengeosys}, porousMultiphaseFoam \cite{horgue2015open}, GeoChemFoam (focused on pore-scale) \cite{maes2021geochemfoam},  DuMux \cite{flemisch2011dumux}, MRST \cite{lie2019introduction}. 
Our work differs from these packages in a number of ways and provides unique contributions to the community, such as:
    \begin{enumerate}
    \item fully integrated workflow from geostatistics to complex flows, including post- and co-processing;
    \item focus on a clear and simple mathematical formulation and general-purpose PDE solvers with extra features developed into object-oriented external modules;
    \item novel numerical methods for equation coupling and stable finite volume formulation;
    \item pluri-Gaussian discontinuous fields, multi-scale features, and dual-porosity models;
    \item Out-of-the-box parallel scalability and \texttt{C++} implementation.
	\end{enumerate}

The paper is organised as follows, Section~\ref{sec:models} presents the mathematical formulation of flow and transport models in porous media at the Darcy scale. Section~\ref{sec:RFG} describes the generation of thresholded Gaussian random fields. Then the numerical details of the implementation are  discussed in Section~\ref{numerics} and illustrated with some examples in Section~\ref{sec:examples}. Finally, we present some conclusions.

	\section{Flow and transport models}
	\label{sec:models}

	\subsection{Darcy flow}
In all our models we assume the validity of the Darcy's law for porous medium fully saturated by an incompressible fluid, which reads as follows:
\begin{align}
  \label{eq:Darcy}
   \vel = -\frac{\perm}{\mu} \left(\nabla p - \rho \gb \right), 
\end{align}
where $\vel$ $[LT^{-1}]$ is Darcy velocity, $\perm=\perm(\xb)$ $[L^{2}]$ is the permeability tensor, $\mu$ $[ML^{-1}T^{-1}]$ is fluid viscosity, $p$ $[ML^{-1}T^{-2}]$ is pressure, $\rho$ $[ML^{-3}]$ is the fluid density, $\gb$ $[LT^{-2}]$ is the acceleration of gravity vector.
We assume that $\rho$ and $\mu$ are non constant but dependent only of the solute concentration $c$ $[ML^{-3}]$ although extensions to double diffusive models (where there is an additional dependency on temperature), compressible, and multiphase descriptions are currently being developed in the code.

Darcy's law can be rewritten in terms of a reduced pressure $\prgh=p-\rho\gb\cdot\xb$ removing the hydrostatic pressure contribution\footnote{It is worth noting that we use here the (non-constant) density field $\rho$ for the reduced pressure instead of a reference density of the acqueous phase as more commonly done.} as follows:
\begin{align}
  \label{eq:Darcy_prgh}
   \vel = -\frac{\perm}{\mu} (\nabla \prgh + \del{\gb\cdot\xb} \nabla \rho ).
\end{align}
%

	\subsection{Continuity equation}

The continuity equation is the balance of mass of fluid per volume of porous medium. It can be formulated as
\begin{align}
  \label{eq:continuity}
  \frac{\partial \rho\por}{\partial t} + \diverg (\rho \vel) = \rho^{\ast}Q,
\end{align}
where $\phi$ $[1]$ is the medium porosity and there is a sink/source flow $Q$ $[L^{-3}T^{-1}]$ of density $\rho^{\ast}$.

We can expand the first term of \eqref{eq:continuity} as
\begin{align}
  \label{eq:continuity2}
  \frac{\partial \rho\por}{\partial t} = S_{0}\rho \frac{\partial p}{\partial t} +  \por \frac{\partial \rho}{\partial c}\frac{\partial c}{\partial t},
\end{align}
where $S_{0}= \partial \phi/\partial p $ $[M^{-1}L^{2}T]$ is the storativity which takes into account the linearised porous matrix compressibility.

	\subsection{Transport equation}
The transport equation gives the solute mass balance per volume of porous medium. We write the transport in terms of the solute concentration $c$ $[ML^{-3}]$ as
\begin{align}
  \label{eq:transport}
  \frac{\partial \por c}{\partial t} + \nabla \cdot (\vel c) - \nabla \cdot (\phi \mathbf{D} \nabla c) = Q c^{\ast},
\end{align}
where $c^{\ast}$ is the concentration of the sink/source and $\mathbf{D}$ $[L^{2} T^{-1}]$ is the hydrodynamic dispersion tensor
\begin{align}
  \label{eq:DispTensor}
  \mathbf{D} = D_{m} \mathbf{I} + \alpha_{T} ||\vel||\mathbf{I} + (\alpha_{L} - \alpha_{T}) \frac{\vel\otimes\vel}{||\vel||},
\end{align}
with $\mathbf{I}$ being the identity matrix, $D_{m}$ $[L^{2}T^{-1}]$ the diffusion coefficient and $\alpha_{L}, \alpha_{T}$ $[L]$ the longitudinal and transverse dispersion coefficients respectively.

	\subsection{Density and viscosity models}

Fluid viscosity and density can be written as a function of the concentration of dissolved solutes as well as temperature. In the following we will focus on linear relations for  $\rho(c)$ and $\mu(c)$, namely:
\begin{equation}
\label{eq:linear_laws}
f(c) = f_{0} + f' c \,,
\end{equation}
where $f$ is either density or viscosity. A more extensive choice of non-linear models is available in the code (see \cref{app:density_viscosity}). The same structure can be used for temperature-dependent viscosity and densities and extended to simultaneous dependence on concentration and temperature (double diffusive model).

\subsection{Boussinesq approximation}
It is sometimes convenient to rewrite the continuity equation as a condition for the divergence of the velocity field to make use of incompressible flow splitting algorithms.  In this case, \cref{eq:continuity} can be reformulated as:
\be{divu}
S_{0} \frac{\partial p}{\partial t} + \diverg\vel =  -\frac{\rho'}{\rho}\del{\phi\dpd{c}{t}+\vel\cdot\grad{c}} + \frac{\rho^*Q}{\rho},
\ee
where it can be noticed that the right hand side contains the Lagrangian derivative of $c$. Therefore, by either assuming $\frac{\rho'}{\rho}\ll1$ or advection-dominated  solute transport, we obtain the Boussinesq approximation, which, for non-deformable porous skeleton (i.e., $S_{0} = 0$) and in absence of sources/sinks, reduces to the divergence free condition
\be{boussinesq}
\diverg{\vel}=0.
\ee

%

	\subsection{Dual-porosity and multi-continuum models}
The equations above are no longer a good approximation when the porous medium is highly-heterogeneous. For example, the pore space could be composed of large highly permeable pores (such as fractures) connected to a system of narrow low permeability (micro)-pores. In these cases, two separate continuity and momentum equations can be considered for two overlapping continua representing the porous spaces. Under the assumptions of \cref{eq:boussinesq} and \cref{eq:continuity}, the system of equations for the two pressures is:
\be{dual_porosity}
\diverg\left[-\rho \frac{\perm}{\mu} (\nabla \prgh + \del{\gb\cdot\xb} \nabla \rho )\right]=\tau(\prgh,\widehat{\prgh})
\ee

\begin{equation}
\diverg\left[-\widehat{\rho}\frac{\widehat{\perm}}{\mu} (\nabla \widehat{\prgh} + \del{\gb\cdot\xb} \nabla \widehat{\rho} )\right]=-\tau(\prgh,\widehat{\prgh}),
\end{equation}
where $\widehat{\cdot}$ denotes the variables in the second continuum (e.g. fractures). The transfer term between the two continua is in general an unclosed integral term depending on the connectivity between the domains, but it is commonly approximated as a linear transfer term
\be{transfer}
\tau(\prgh,\widehat{\prgh})=-\tau_0 (\prgh-\widehat{\prgh}),
\ee
where $\tau_0$ is a linear transfer coefficient between matrix and fracture continua.

	\section{Random fields generation}
	\label{sec:RFG}


Continuous and thresholded Gaussian random fields (also known as pluri-Gaussian) \cite{hesse2014generating,beucher2016truncated, journel1978mining}  can be generated using different approaches (sequential Gaussian simulations \citep{dimitrakopoulos2004generalized}, Markov chain probability \cite{kemeny1976markov}, multiple-point statistics \cite{strebelle2002conditional}). In contrast with the number of geostatistics open-source toolboxes available (GSLib \cite{deutsch1998gslib}, T-PROGS \cite{carle1999t}), we use an integrated framework to combine geostatistics with flow and transport simulation and post-processing. Another key difference is that our random field generators do not make use of the fast Fourier transform or other discrete transform. This can make the generation more computationally intensive but this is mostly overcome by the efficient \texttt{C++} implementation and parallel scalability. In our tests, in fact, although the random field generation can be significantly expensive, it is nevertheless often negligible compared to the cost of the flow and transport solvers, and these steps can be fully integrated in a single run. 

\subsection{Continuous and pluri-Gaussian truncated fields} \label{sec:pTGS}
\begin{figure}[htbp]
	\centering
	\includegraphics[width=0.9\textwidth]{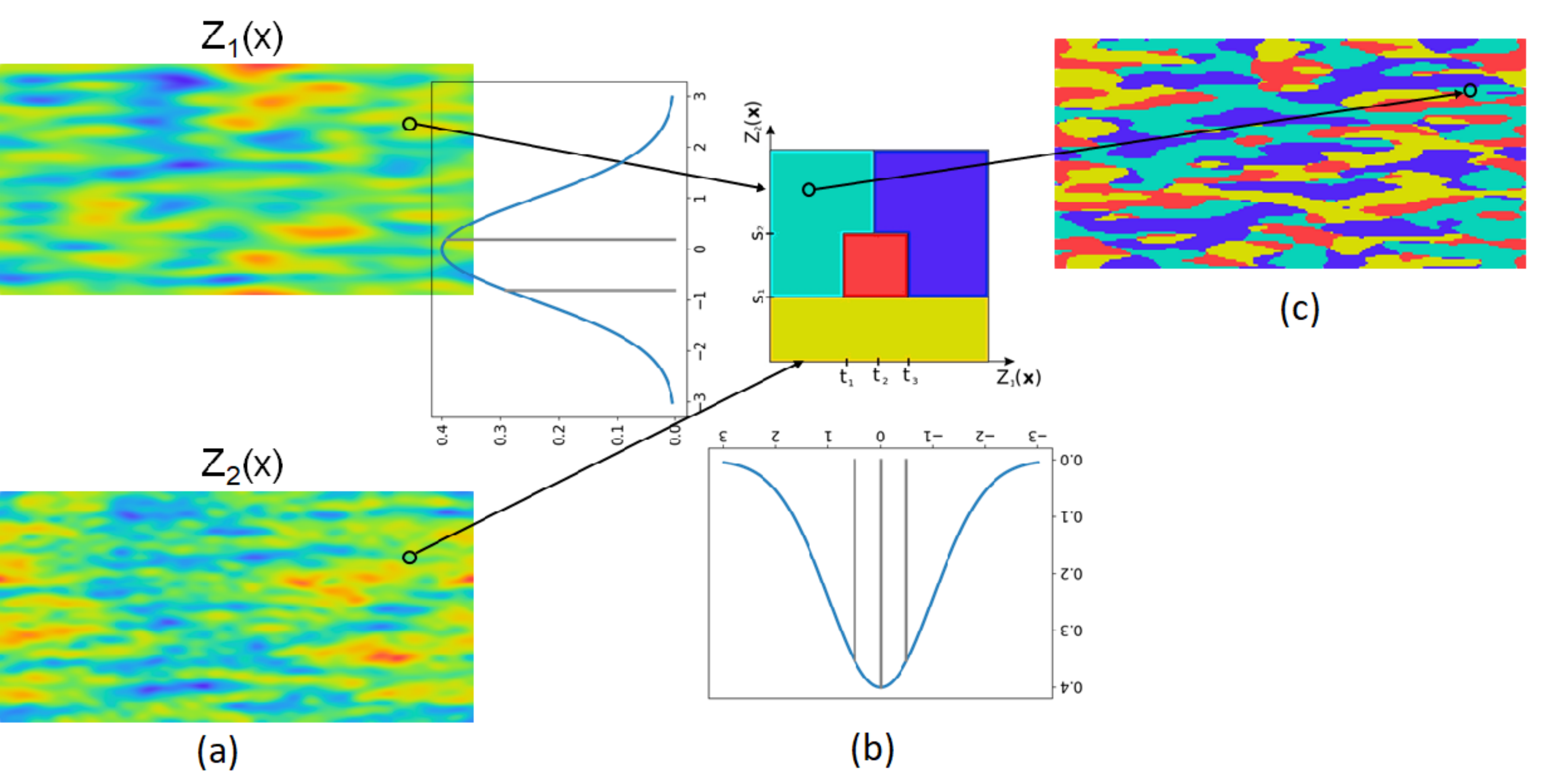}
	\caption{Pluri-Guassian simulation method: (a) two continuous Gaussian random fields are created; (b) at every position in space, the values of the continuous GRFs are used as coordinates to enter the truncation rule; (c) an heterogeneous non-Gaussian field is generated.}
	\label{fig:PGSmethod}
\end{figure}

Gaussian random fields (GRF) have often been adopted in geostatistics to mimic the spatial distributions of geological properties due to their mathematical properties.
Generally these fields well fit the purpose of geological descriptors as long as the spatial transition of the geological properties is smooth. The choice of correlation function can tune the smoothness of the function but always in a continuous manner.
However, sediments' deposits rarely show a continuous structure. In these cases a continuous GRF can be post-processed with thresholding or binning to allow for abrupt transitions and discontinuities. This process is called single-truncation rule.
However, if the sediment pattern is characterised by non-layered or non-stratified geometries where one category share its border with more than other two, the single-truncation model does no longer provide realistic results. To reproduce non-stratified geometries the thresholding needs to be done on a multivariate Gaussian random field. When two (or more) independent Gaussian random fields are thresholded according to a multi-dimensional truncation rule, the resulting field is called a pluri-Gaussian random field. The underlying idea is to simulate two or more GRFs on a domain and compare them through a series of inequalities which allows us to assign a unique value to each cell (Fig.~\ref{fig:PGSmethod}).

\begin{figure}[htbp]
	\centering
	\includegraphics[width=0.3\textwidth]{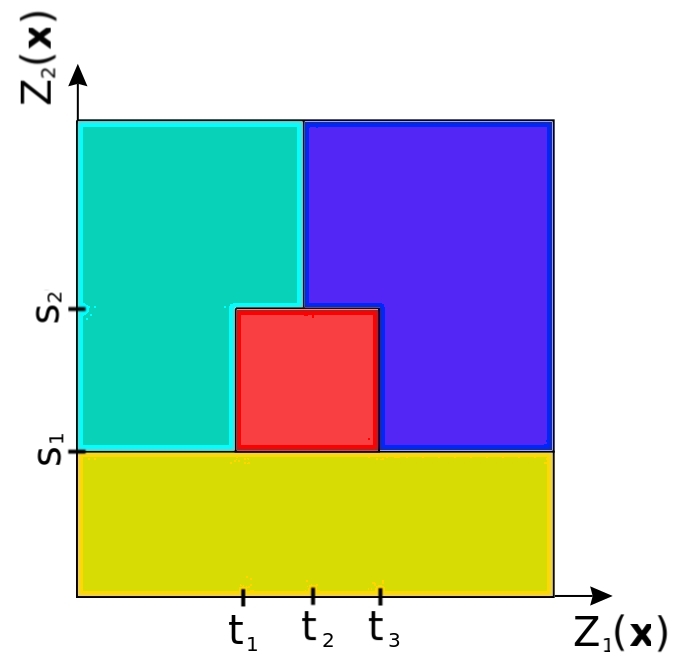}
	\caption{Qualitative truncation rule. Note that the areas of the different facies do not correspond to the their respective proportions in the simulations, not because of the qualitative nature of the image, but because the underlying continuous variables are not uniformly distributed but Gaussian.}
	\label{fig:TruncRule}
\end{figure}

An example of a two-dimensional \textit{truncation rule} with three ($t_1,t_2,t_3$) and two ($s_1,s_2$) thresholds in the fist and second dimension, respectively, applied to two independent random realisations of GRF ( $Z_1(\boldsymbol{x})$ and $Z_2(\boldsymbol{x})$) is depicted in Fig.~\ref{fig:TruncRule}. In the example, all red cells in the domain satisfy the condition $t_1<Z_1(\boldsymbol{x})<t_3$ and $s_1<Z_2(\boldsymbol{x})<s_2$. The rule gives rise to four facies (geological domains) all connected to each other.

If the thresholds are expressed in terms of percentiles of the GRF, the probability, i.e. the proportion of the facies $i,j$ is
\begin{equation}
	p_{i,j} = \left[ G_1(t_i) - G_1(t_{i-1}) \right] \left[ G_2(s_j) - G_2(s_{j-1}) \right] .
\end{equation}
Given $I$ and $J$ thresholds, respectively, the problem of choosing a truncation rule for $N$ facies is underdetermined if $N<(I+1)(J+1)$. Extra constraints might come from the connectivity or surface area of each facies.
In the previous examples a total of twelve regions are identified by the Cartesian truncation diagram but these have been then merged into a total of four facies.
This can be solved by constrained minimisation problem involving $(I+1)(J+1)$ error functions. In other words, in this case, there are infinite threshold combinations which honour a given set of volume proportions and some constraints need to be defined beforehand to allow for a single solution. A more detailed analysis of this problem can be found elsewhere \cite{mariethoz2009truncated}.
 In the computational framework we present here, the user needs to specify directly the thresholds rather than the facies volumes.

\begin{sidewaysfigure}[htbp]
\centering
  \includegraphics[width=0.25\linewidth]{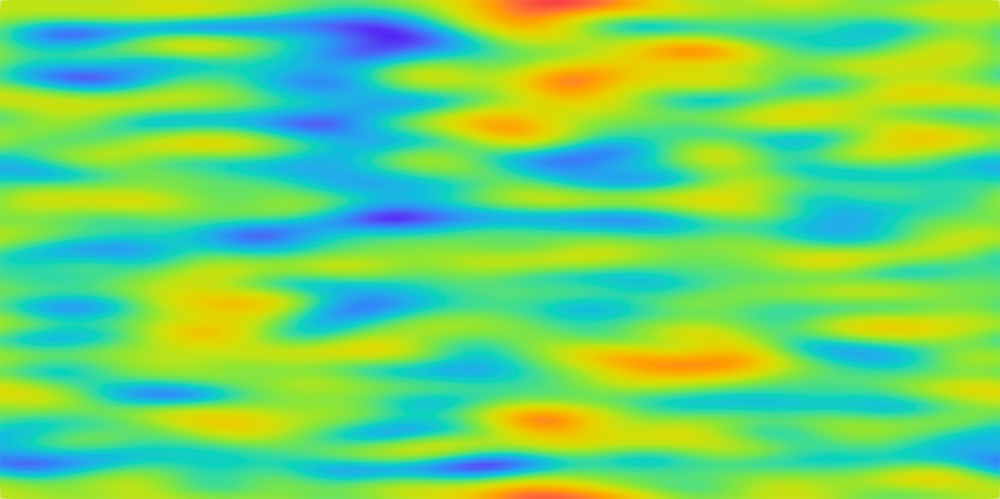}
  \hspace{0.5cm}
  \vspace{0.2cm}
  \includegraphics[width=0.25\linewidth]{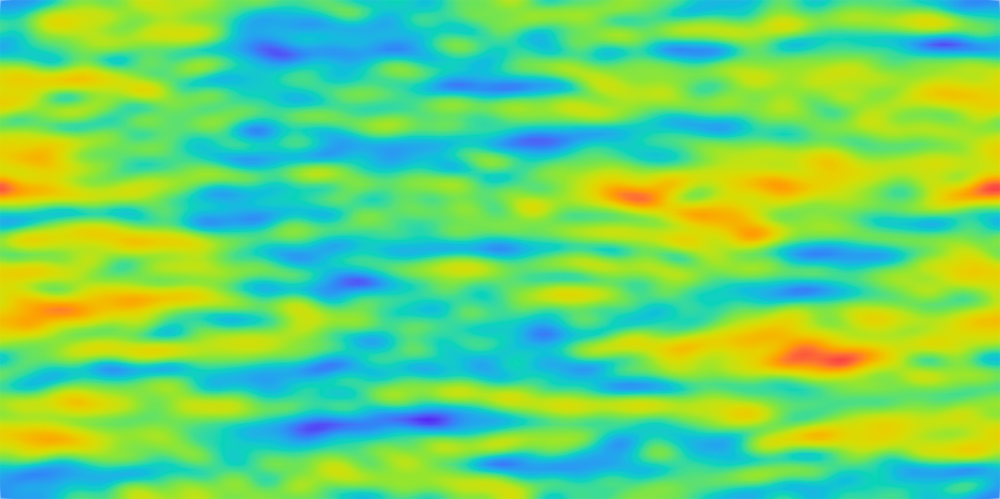}
  \hspace{0.5cm}
  \vspace{0.2cm}
  \includegraphics[width=0.25\linewidth]{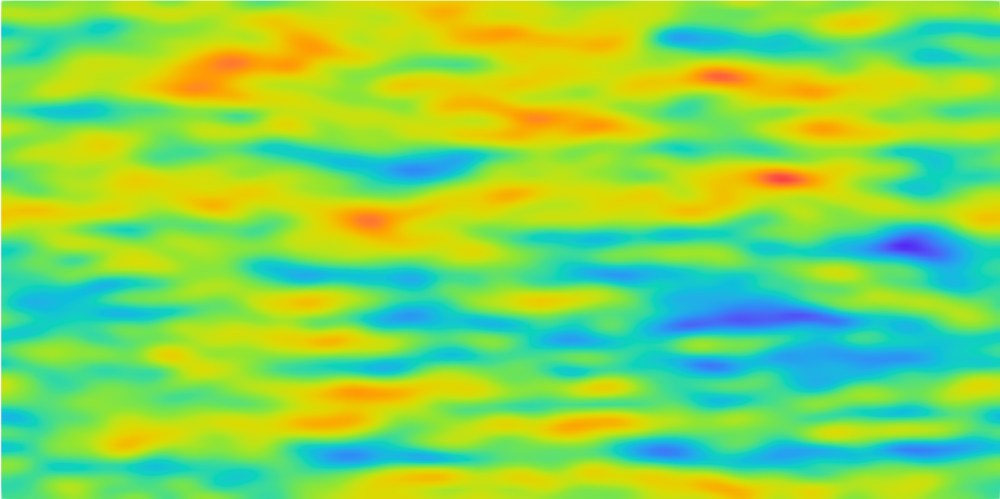}
  \vspace{0.2cm}
  \includegraphics[width=0.25\linewidth]{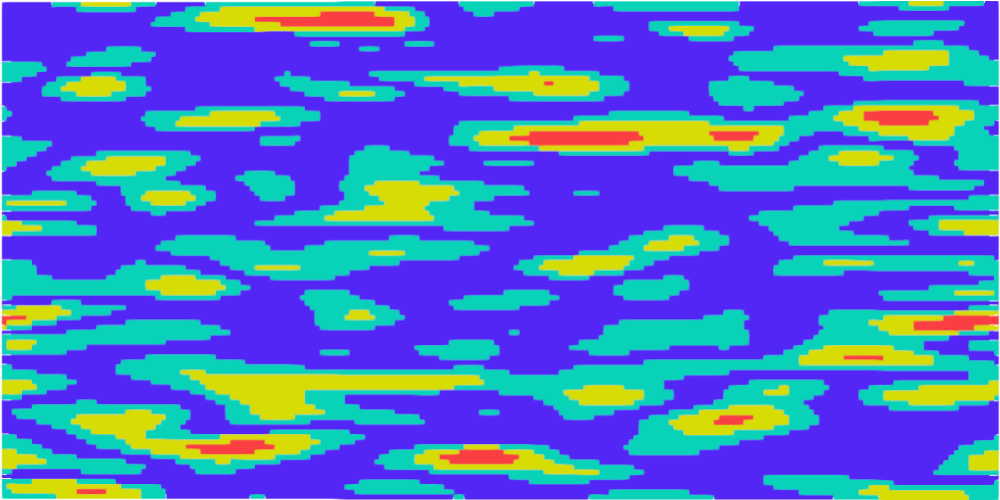}
  \hspace{0.5cm}
  \vspace{0.2cm}
  \includegraphics[width=0.25\linewidth]{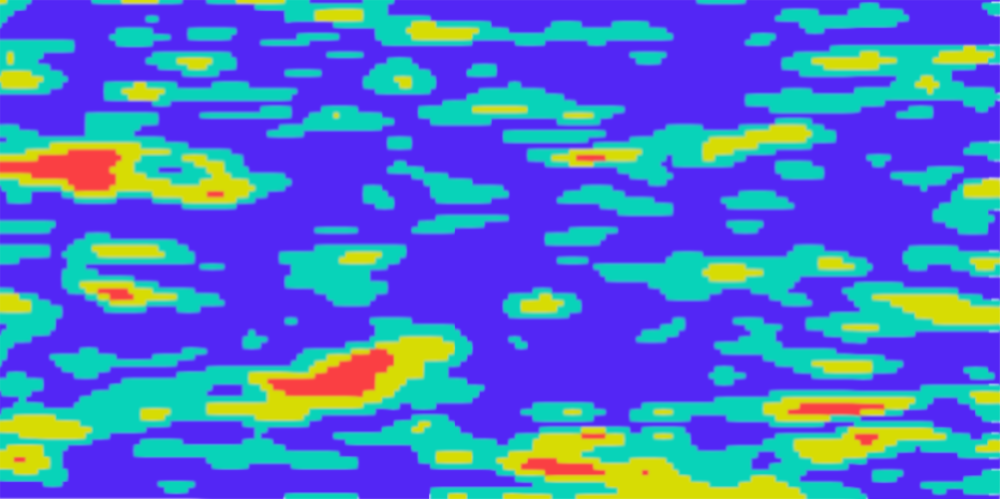}
  \hspace{0.5cm}
  \vspace{0.2cm}
  \includegraphics[width=0.25\linewidth]{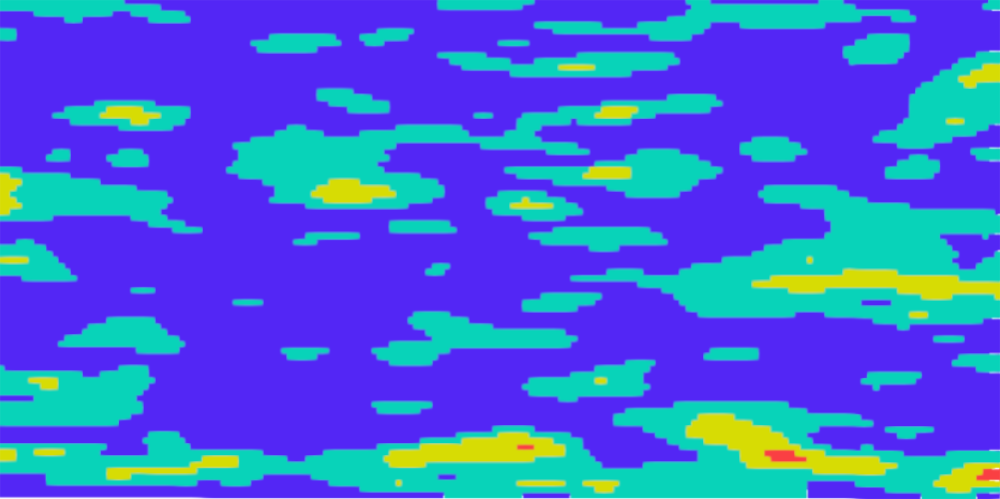}
  \vspace{0.2cm}
  \includegraphics[width=0.25\linewidth]{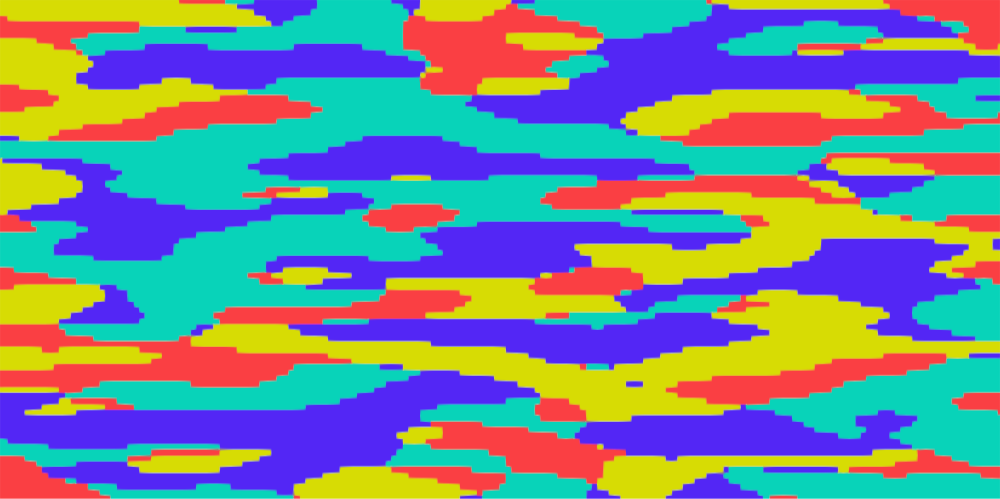}
  \hspace{0.5cm}
  \includegraphics[width=0.25\linewidth]{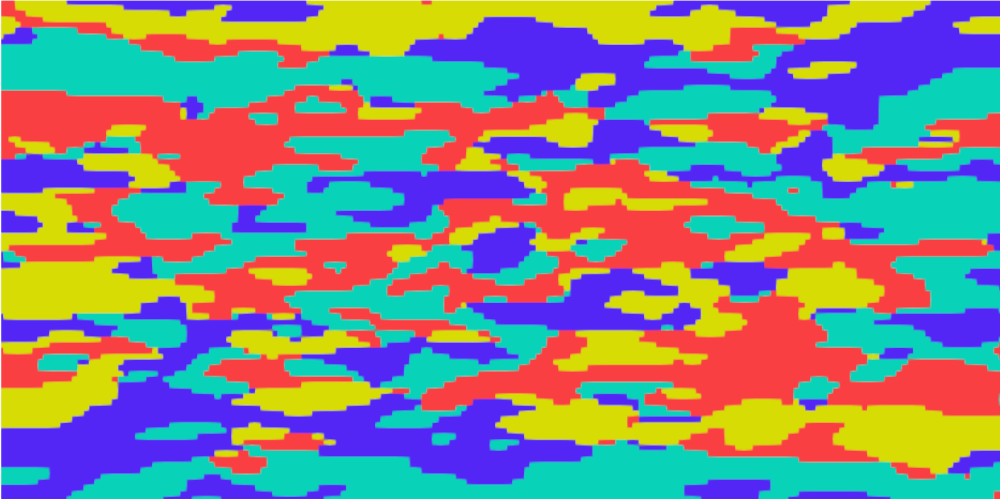}
  \hspace{0.5cm}
  \includegraphics[width=0.25\linewidth]{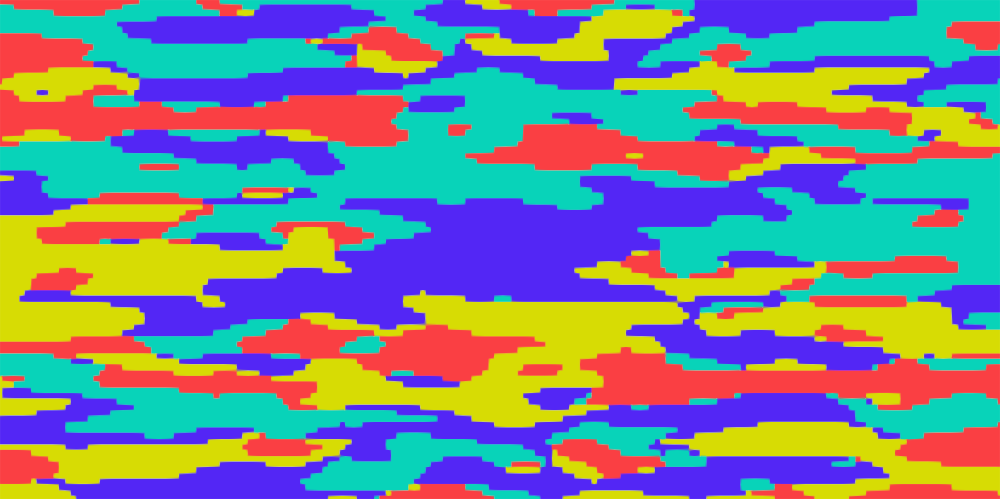}
\caption{ Two-dimensional lognormal random field with (a) Gaussian correlation (b) exponential correlation and (c) Matérn correlation with $\nu=1$ (see appendices for the mathematical definition of these covariances)}
\label{fig:contCov}
\end{sidewaysfigure}

As an example, in Figure \ref{fig:contCov}, we show two-dimensional realisations of continuous, truncated and pluri-truncated GRFs obtained using formula \eqref{eq:DiscreteFT} and Gaussian \eqref{eq:gauSpe}, exponential \eqref{eq:expSpe} and Matérn \eqref{eq:matSpe} covariances.

	\section{Numerical implementation}
	\label{numerics}
In this section we present a few more details about some of the most important solvers and tools implemented in \texttt{SECUReFoam}, namely:
\begin{itemize}
\item \textsf{setRandomField}: a pre-processing utility to generate random fields for permeability, porosity and other properties;
\item \textsf{rhoDarcyFoam}: a variable-density flow and transport solver;
\item \textsf{dualSimpleDarcyFoam}: a dual-porosity Darcy solver;
\end{itemize}
together with boundary conditions, sources and post-processing tools. Other solvers and utilities included the library but not described here include
\begin{itemize}
\item \textsf{simpleDarcyFoam}: a simple single-phase Darcy solver;
\item \textsf{multiRateScalarTransportFoam}: a multi-rate mass transfer model \cite{municchi2020generalized,municchi2021heterogeneous};
\item \textsf{poroelasticFoam}: a linear Biot poroelasticity solver;
\item \textsf{dualRhoDarcyFoam}: a dual-porosity variable density solver;
\item mesh importing utilities from \textsf{gslib, ijk,} and \textsf{grdecl} format  ;
\item \textsf{spatialPdf}: an utility to create linear and log-spaced histograms from spatial data;
\item \textsf{fieldMetrics}: co- and post-processing utility to compute statistics of spatial data. This is a specialised and extended version of the \textsf{functionObject} structure of \of.
\end{itemize}


	\subsection{Meshing and random field generator \textsf{setRandomField}}

\of\ is an unstructured finite-volume library that can deal with arbitrary cell types. Cartesian grids are treated internally as non-structured meshes. In the present work we have limited ourselves to Cartesian meshes generated by the \of\ native \textsf{blockMesh} utility. However, included in \texttt{SECUReFoam} are utilities for importing corner-point reservoir models, based on the work done within the MRST project \cite{lie2019introduction}. Alternatively more complex geometries can be meshed combining from elementary volume objects or CAD files.

Once the grid has been created, the material properties can be populated with the random field generator \textsf{setRandomField}. As we have seen in the previous section, this is based on the generation of GRFs.
Several well-known methods have been implemented in literature for generating GRFs \cite{hesse2014generating, zhang2006random, deutsch1998gslib, ruan1998efficient}. According to Mandelbrot and Van Ness (1968) \cite{mandelbrot1968fractional}, a GRF can be represented using a stochastic Fourier integral:
\begin{equation} \label{eq:stoFouInt}
	Z(\mathbf{x}) = \int_{-\infty}^{+\infty} e^{-2 \pi i \mathbf{a} \cdot \mathbf{x}} \sqrt{S(\mathbf{a})} dW(\mathbf{a}),
\end{equation}
where $\mathbf{a}$ are frequencies, $dW(\mathbf{a})$ is a complex-valued white noise random measure and $S(\mathbf{a})$ is the amplitude of the spectral measure.
The latter can be found as the Fourier spectrum of correlation functions written in spherical coordinates. To allow for deformations along the $x,y,z$ directions (i.e., different correlation lengths in each direction), a rescaling is implemented to the three-dimensional frequency space. This can be easily extended by including rotations, or more general transformation matrices, allowing therefore more flexibility in the orientation of the facies. 

The discrete form of \eqref{eq:stoFouInt} becomes
\begin{equation} \label{eq:DiscreteFT}
	Z(\mathbf{x}) = \sum_{j=0}^{N_f} \cos(2 \pi \mathbf{a}_j\cdot \mathbf{x}) \sqrt{S(\mathbf{a}_j)} W_j + i 
\sum_{j=0}^{N_f} \sin(2 \pi \mathbf{a}_j \cdot\mathbf{x}) \sqrt{S(\mathbf{a}_j)} W'_j.
\end{equation}
Instead of relying on the FFT algorithm to compute \eqref{eq:DiscreteFT}, which would significantly reduce the computational cost of computing the Fourier integral, we have directly implemented the formula for a number of reasons. First of all, FFT can hardly be extended to non-Cartesian unstructured meshes. These are very common in geology and reservoir simulation. The other disadvantage of the FFT is that the random realisation becomes dependent on the spatial discretisation used. This makes (deterministic) convergence studies (for a given random field) hardly achievable. With \eqref{eq:DiscreteFT} instead, for a given (pseudo-)random set of Gaussian variables $W_i$, we can generate multiple random fields at different spatial resolution. Finally this approach allows us to relax the periodicity assumptions by including larger wave numbers to the sum, simulating therefore non-stationary random fields.

Once the GRFs are generated as above, they are either scaled to obtain the desired mean and variance, or thresholded, as explained in \cref{sec:pTGS}, to obtain discontinuous fields. The options and parameters for the random field generation are described in \cref{app:setRandomFieldDict}.

	\subsection{Variable density solver \textsf{rhoDarcyFoam}}
	\label{sec:rhodarcy}

\textsf{rhoDarcyFoam} solves variable-density flow and transport problems with the models presented in \cref{sec:models}. Different formulation of the equations can be used, based on standard pressure \eqref{eq:Darcy} or the reduced pressure \eqref{eq:Darcy_prgh}, with or without Boussinesq approximation (\eqref{eq:boussinesq} vs \eqref{eq:divu}). The transport model is based on \eqref{eq:transport}. Multiple parameters can be defined as heterogeneous fields, including permeability, porosity, dispersion parameters, and storativity.
Various dispersion, density, and viscosity models are included by means of a object-oriented modular structure and selectable through simple input file (following \of\ standard, most of the physics settings have been included into the \textsf{transportProperties} dictionary).

The resulting system is a non-linear coupled PDE system. An under-relaxed Picard iteration loop is employed to couple the transport and Darcy equations. This is based on the \textsf{pimple} class in \of, that allows the control of outer iterations to exit the loop when the residuals fall below a certain value or when the maximum number of iterations is reached. This also allows to remove the relaxation in the last iteration to prevent it from delaying the dynamic evolution of the fields.
Using existing \of\ capabilities, fully implicit (through Picard iterations) and explicit time stepping (by setting one single iteration) can be implemented with first and second-order backward integration. Adaptive time-stepping is included based on a maximum Courant number computed using the total fluxes. The whole formulation is based on the fluxes across faces, rather than the cell-averaged velocities. This ensures exact mass conservation. Optionally, each linear system can be solved multiple times to better incorporate explicit corrections due to non-orthogonal meshes or flux-limiter schemes, as is the case for most of the standard solvers shipped with \of.

The permeability and dispersivity field can be specified as a symmetric matrix field on the faces (e.g. inverse of face transmissibility) or at the faces. While the former does not require any interpolation, the latter instead requires an interpolation to the faces. The standard interpolation schemes in \of\ are overridden to enforce weighted harmonic interpolation for diffusive fluxes. Gradients at the faces are approximated with two-point or multi-point flux approximation (TPFA, MPFA). The latter is based on least-squares high-order gradient  approximation, included in \of.

	\subsection{Dual-porosity solver \textsf{dualSimpleDarcyFoam}}
Dual-porosity solvers are available for simple constant density Darcy flow (\textsf{dualSimpleDarcyFoam}) as well as an extension of the variable density solver for dual-porosity problems (\textsf{dualRhoDarcyFoam}). We limit ourselves here to the description of the Darcy solver.
The dual-porosity model \eqref{eq:dual_porosity} is solved in a segregated (iterative) manner. To increase the convergence of the inner iteration, we make use of a new splitting scheme \citep{splitting} recently developed for coupled systems. These are based on an approximate Schur-complement that allows to increase the coupling between the equations.

	\subsection{Boundary conditions and other modules}
Various new elements as separate modules to be used with various solvers and applications. This includes many new boundary conditions, source terms, and post-processing utilities.
Most new boundary conditions have been implemented based on a new general-purpose Robin boundary condition (\textsf{Robin}) with linear \cite{boccardo2018robust} or exponential \cite{pimenta2019coupled} reconstruction of the solution near the wall. This is then extended to include advective fluxes (\textsf{RobinPhi}) and then used as base class for more complex BCs.  All the BCs are implemented using the standard \of\ class structure for \textsf{fvPatchFields}. We describe here two particular types of BCs which will be extensively used in the numerical results.

	\paragraph{Darcy-based boundary conditions.}
	Specific boundary conditions have been implemented for single-phase and variable-density Darcy solvers. These include a \textsf{darcyFixedVelocity} condition to impose a pressure gradient to ensure a fixed inlet or outlet velocity is obtained. This automatically switches between different formulations of the pressure equation, anisotropic permeabilities, different solvers, and the presence of gravity. \textsf{hydrostaticPressure} is a condition that allows to specify the dynamic pressure while using a solver that includes the hydrostatic pressure in the total pressure.
		
	\paragraph{Flux-based boundary conditions.}
In practical applications it is often convenient to specify a boundary condition for the total fluxes at a boundary for flow and transport. To this aim, we developed the boundary conditions named \textsf{darcyFixedVelocity} and 
 \textsf{fixedTransportFlux}. The former adapts automatically the specific formulation used by the solver (see \cref{sec:rhodarcy}) to impose a given fixed velocity. The latter
  can be written in general, for a given field $f$,
\begin{equation} \label{eq:fixedTransportFlux}
 -D \nabla_n f  + u_n f  = u_n f_a - d (f_d-f)/\Delta + F,
\end{equation}
where, the first and second terms of the left hand side represent the diffusive and advective normal fluxes, respectively, and the left hand side contains the external fluxes imposed by the user. These are an advective flux with a given inlet concentration $f_a$ and the true fluid velocity $u_n$, a discrete diffusive flux computed with a diffusion coefficient $d$ and an external concentration $f_d$ (that accounts for the given half-cell distance from the boundary $\Delta$), and an explicit flux $F$.

	\section{Numerical examples}
	\label{sec:examples}
        In the following we present a series of numerical example to illustrate the capabilities of \texttt{SECUReFoam}. The heterogeneous and dispersive Henry problem (subsection~\ref{subsec:Henry}) shows how the platform deals with variable density in heterogeneous porous media. The performance of the under-relaxed Picard is demonstrated by solving a strongly non-linear viscous fingering instability (subsection~\ref{subsec:viscfing}) and a highly heterogeneous Horton-Rogers-Lapwood problem (subsection~\ref{subsec:HRL}), which also shows the use of truncated pluri-Gaussian fields. Finally, we present a quarter five-spot injection problem to demonstrate the dual-porosity formulation and the capability to deal with discontinuous highly heterogeneous fields.

	\subsection{Heterogeneous Henry problem}\label{subsec:Henry}
	The Henry problem \citep{Henry1964} is an abstraction of the seawater intrusion in a coastal aquifer. It has been extensively used to understand the interaction between the aquifer and the sea and as benchmark of variable density groundwater flow codes \citep{Simpson2004}. In the Henry problem the aquifer is represented by a $2 \times 1$~m rectangle. The right side of the domain is the boundary with the sea where hydrostatic pressure is prescribed using seawater density for the flow equation and an inlet/outlet boundary condition is used for the transport equation. On the left boundary a freshwater flow is simulated by prescribing flow $Q_{\text{in}}$. The rest of the boundaries are impervious. Density depends linearly on the salt concentration $\rho(c) = \rho_{0} + \beta c$ and viscosity is constant. The boundary conditions for the concentration are homogeneous Neumann everywhere but the right (sea) boundary where a constant Dirichlet ($c=1$) is imposed. However, to avoid the formation of a boundary layer when an outward flux develops, an automatic switching is adapted to impose Neumann condition when the velocity is pointing outwards.

        \begin{table}[htb]
          \centering
          \begin{tabular}{lll}
            \hline
            Parameter & Value & Units\\
            \hline
            $Q_\text{in}$ & 6.6e-2 & kg/s\\
            $\phi$ & 0.35 & -\\
            $\rho_{0}$ & 1000 & kg/m$^3$\\
            $\beta$ & 0.6832 & -\\
            $c_{\text{sea}}$ & 36.5925 & kg/m$^3$\\
            $k$ & 1.020408e-9 & m$^2$ \\
            $\mu$ & 0.001 & Pa s\\
            $D_{m}$ & 6.6e-6 & m$^2$/s\\
            $\alpha_{L}$ & 0.1 & m \\
            $\alpha_{T}$ & 0.01 & m \\
            $\lambda_{x}$ & 0.2 & m \\
            $\lambda_{z}$ & 0.1 & m \\
            $\sigma^{2}_{\log k}$ & 2 &\\
            \hline
          \end{tabular}
          \caption{Henry problem parameters for the diffusive and dispersive cases.}\label{tab:Henry}
        \end{table}

        \begin{figure}[!hbtp]
          \centering
          \includegraphics[width=0.49\linewidth]{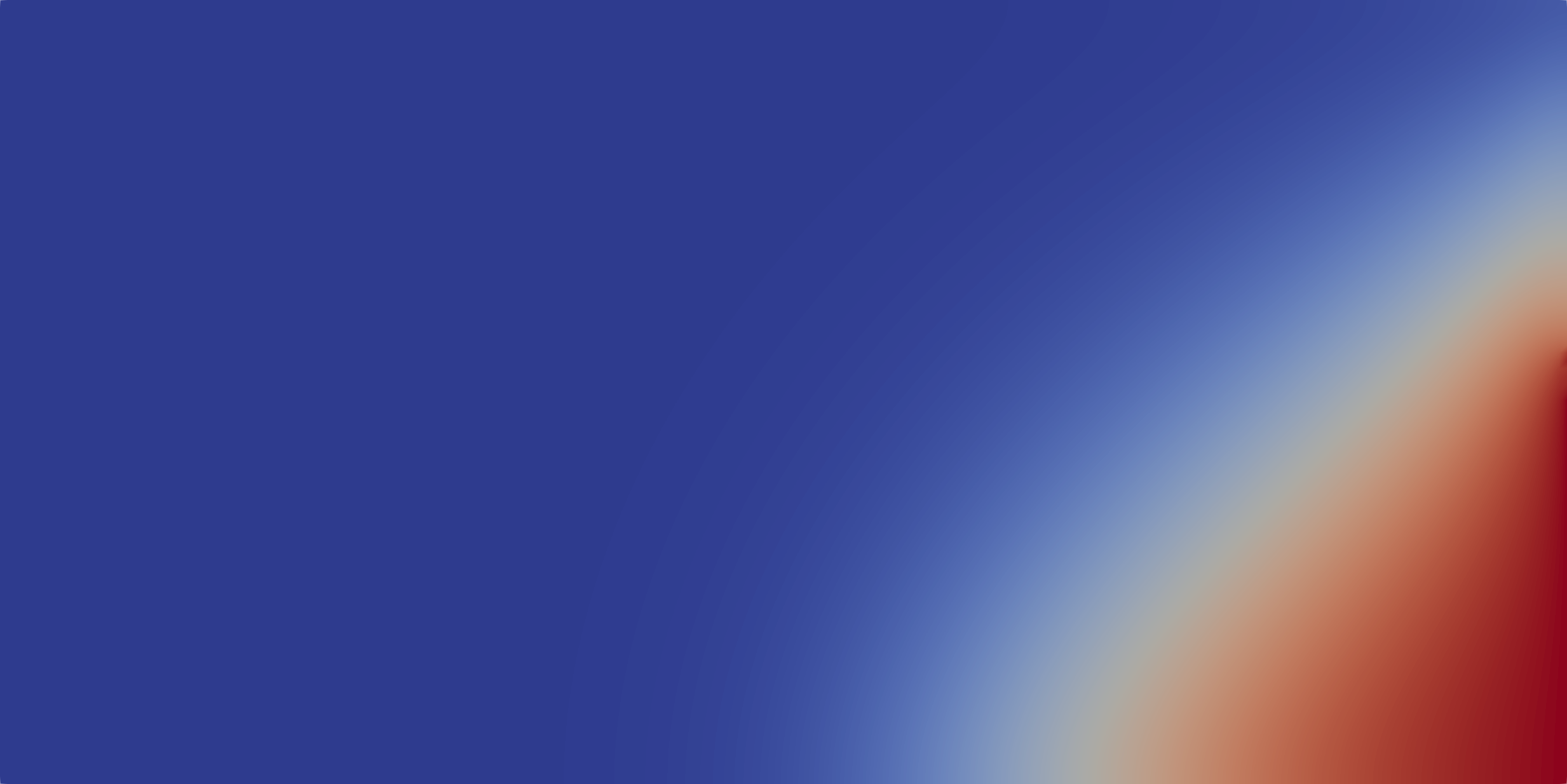}
          \includegraphics[width=0.49\linewidth]{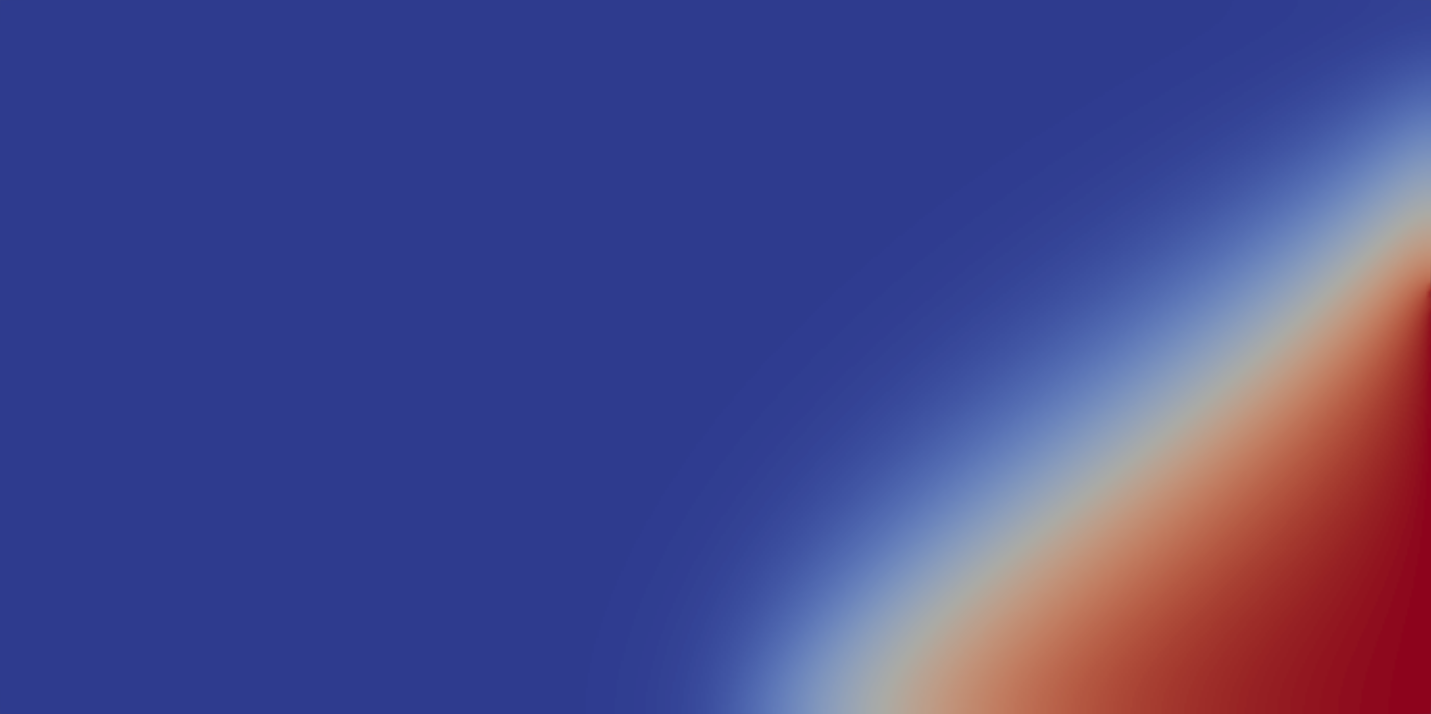}
          \caption{Solution for the concentration for the original Henry (left) and dispersive Henry problem (right) at time 6000 s. As a result of hydrodynamic dispersion the freshwater-saltwater interface is flatter and advances more inland.}\label{fig:Henry}
        \end{figure}

        \begin{figure}[!hbtp]
          \centering
          \includegraphics[width=0.49\linewidth]{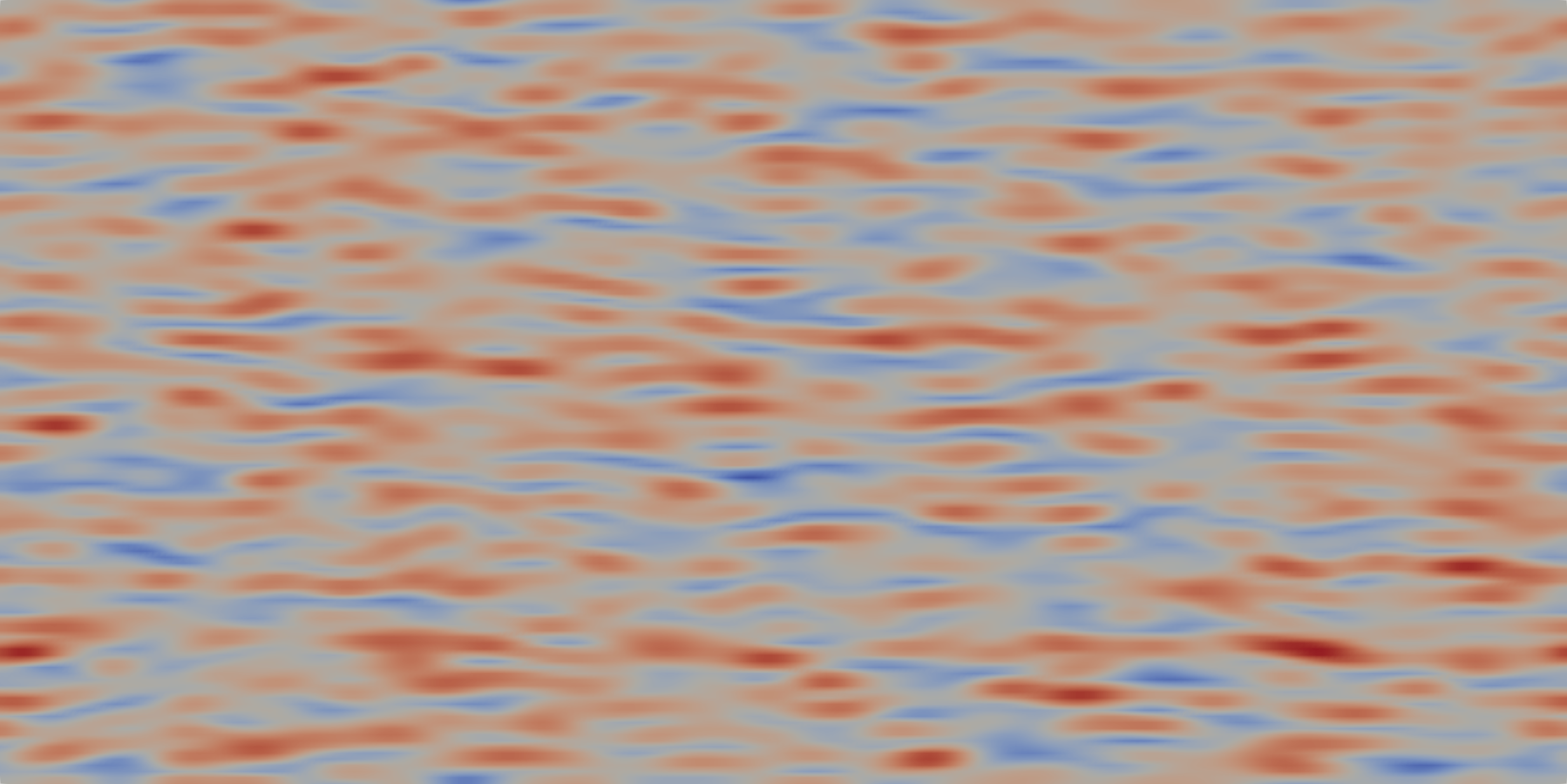}
          \includegraphics[width=0.49\linewidth]{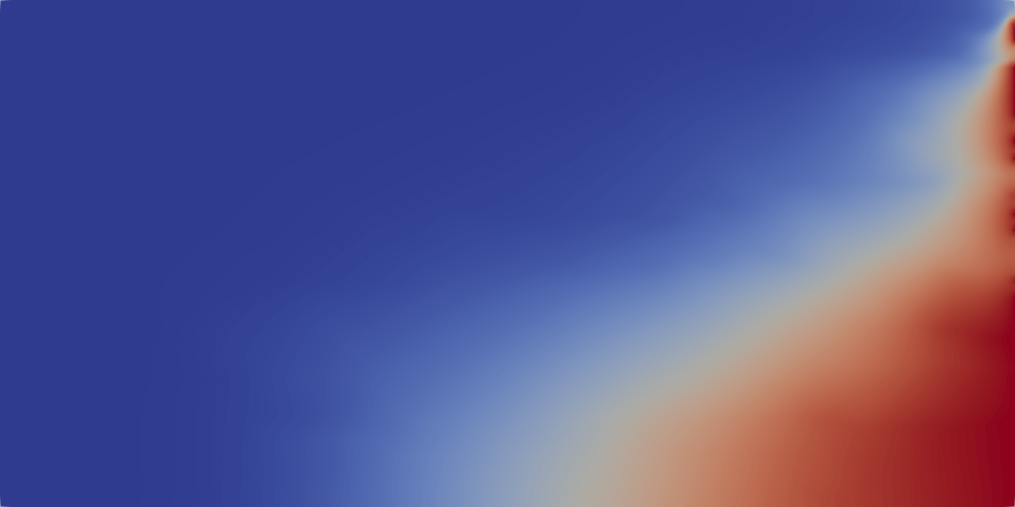}
          \caption{Map of log-permeabiity field (left) used to solve the heterogenous Henry problem permeability and the resulting concentration (right) at time 6000 s.}\label{fig:Henry2}
        \end{figure}

Along with the original Henry problem (figure \ref{fig:Henry} left), we also consider a dispersive case (figure \ref{fig:Henry} right) and a heterogeneous case (figure \ref{fig:Henry2}) to illustrate the role of hydrodynamic dispersion and the heterogeneity of the permeability. The parameters for all cases can be found in Table~\ref{tab:Henry}. As shown in figure \ref{fig:Henry}, hydrodynamic dispersion causes the interface between the seawater and the fresh water to become flatter. Dispersion also affects the movement of the saltwater wedge, which travels further inland in the dispersive case than in the diffusive case. The same effect on the wedge's toe position is observed for the heterogeneous case. Heterogeneity, however, distorts more strongly the geometry of the saltwater-seawater interface. The low permeability zones near the sea (right) boundary modify the discharge of the freshwater and lead to the formation of high concentration zones along the boundary. The boundary conditions on that boundary are switching between Dirichlet and Neumann, therefore creating a seemingly oscillatory patterns which is however stable and purely due to the heterogeneities and the resulting horizontal velocity which switches between positive and negative. This also shows the limitation of this simple testcase in highly heterogeneous systems.

	\subsection{Viscous fingering}\label{subsec:viscfing}

\begin{figure}[!hbtp]
            \includegraphics[width=0.49\linewidth]{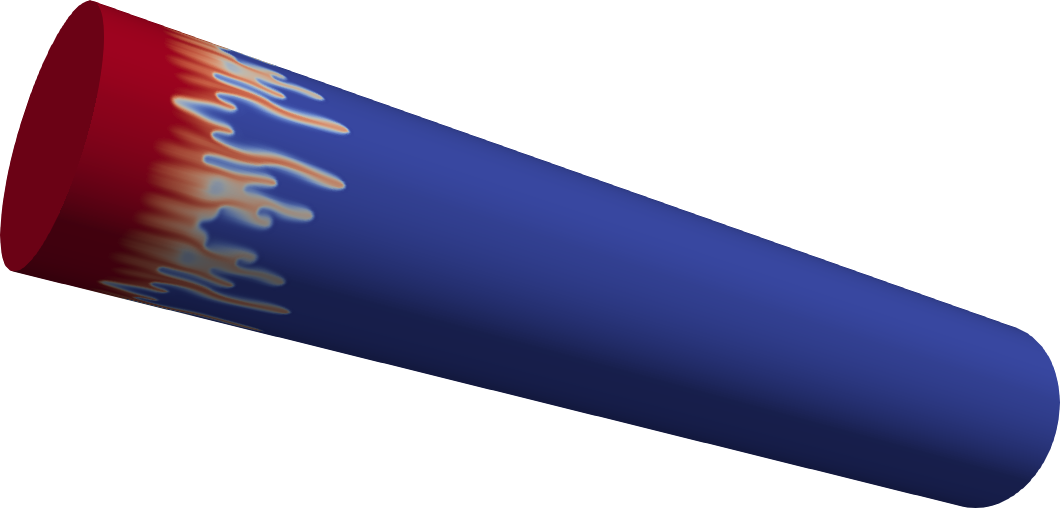}
          \includegraphics[width=0.49\linewidth]{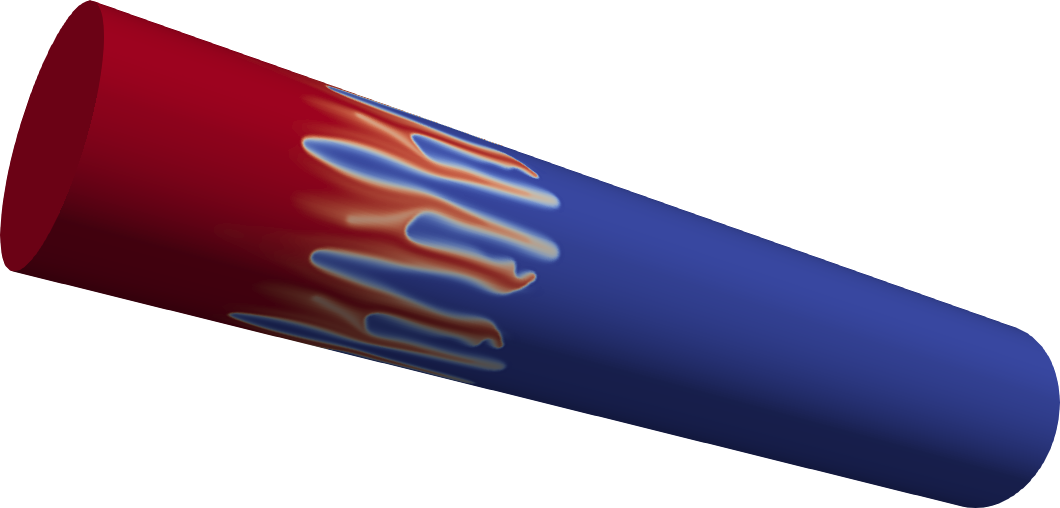}
          \caption{3D viscous fingering  at $Pe= 10^{3}$ at dimensionless times 0.3 and 0.7. A high viscosity fluid (red) displaces a low viscosity fluid forming a fingering instability. As time passes the initially small fingers merge into larger ones. Shielding and tip splitting can be observed in the fingers of the panel on the right.}\label{fig:vf3D}
        \end{figure}

        Viscous fingering is a flow instability that appears when a fluid displaces another one of different viscosity \citep{Homsy1987}. The instability is caused by the difference in mobility between the fluids, which leads to the formation of fingering patterns. The patterns displays a complex tip-splitting, shielding and coalescence dynamics, which is affected by the medium heterogeneity \citep{Araktingi1993, DeWit1997, Nicolaides2015} and chemical reactions \citep{DeWit1999, Rose2013, Kim2021}. Figure~\ref{fig:vf3D} shows an example of miscible viscous fingering in a 3D cylindrical geometry. In this case viscosity is a function of concentration $\mu(c) = \mu_{0} e^{Rc}$ with $R=-3$. A fluid with $c=1$ is injected at constant rate from the left boundary displaces the resident fluid ($c=0$). Pressure and a zero gradient for concentration are prescribed at the outlet. The system is characterised by the Péclet number $Pe =q_{0}L/\phi D_{m} = 10^{3}$, where $q_{0}$ is the velocity at the inlet, $L$ the domain length, $\phi$ the porosity and $D_{m}$ the diffusion coefficient. We observe how the interface deforms and the fingering pattern appears (figure~\ref{fig:vf3D}, left). At late times (figure~\ref{fig:vf3D}, right), the fingers merge and the pattern is form by fewer thicker fingers. The tip splitting and shielding mechanisms are also reproduced in the simulation.

These results have been obtained with $200\times200\times1000$ Cartesian grid which is cut and adapted to the cylinder walls using the \of native meshing tool \textsf{snappyHexMesh}. The adaptive time step is chosen to be proportional to the local mesh size and the inverse of the local velocity magnitude, with a proportionality constant of $0.1$. The instability patterns are however only qualitatively independent on the mesh size and time step as no initial perturbation is imposed. Particularly important here (and even more for heterogeneous permeability fields) is the harmonic interpolation of the viscosity and diffusive coefficients to properly characterise the fingering.

	\subsection{Unstable flow in highly heterogeneous media}\label{subsec:HRL}
        The Horton - Rogers - Lapwood (HRL) problem \citep{Horton1945, Lapwood1948} is a heat transport problem in which a temperature difference is prescribed between the top and bottom boundaries of a rectangular domain. Fluid density decreases linearly with temperature so that an unstable density stratification is form that triggers a Rayleigh - Bénard instability. The system is characterised by the Rayleigh number
        \begin{align}
          \label{eq:Ra}
          Ra = \frac{k  g \Delta \rho H}{\phi \mu K_{th}}
        \end{align}
        where $K_{th}$ is the thermal conductivity, $k$ permeability, $\mu$ viscosity, $\Delta \rho$ the density difference between the top and bottom boundaries, $H$ the height of the domain and $\phi$ the medium porosity. For a square domain the system is stable for $Ra < 4\pi^{2}$. For $ 4\pi^{2} < Ra < 1300$ the system becomes unstable and convection cells that occupy the whole domain form \citep{Graham1994}. For $Ra > 1300$ the convection regimes becomes chaotic and flow organises itself in columnar patterns \citep{Hewitt2013}. During the convection-dominated regimes, mixing and heat fluxes through boundaries are significantly increased with respect to the stable regime \citep{Hidalgo2018}. 

\begin{figure}[!hbtp]
          \includegraphics[width=0.49\linewidth]{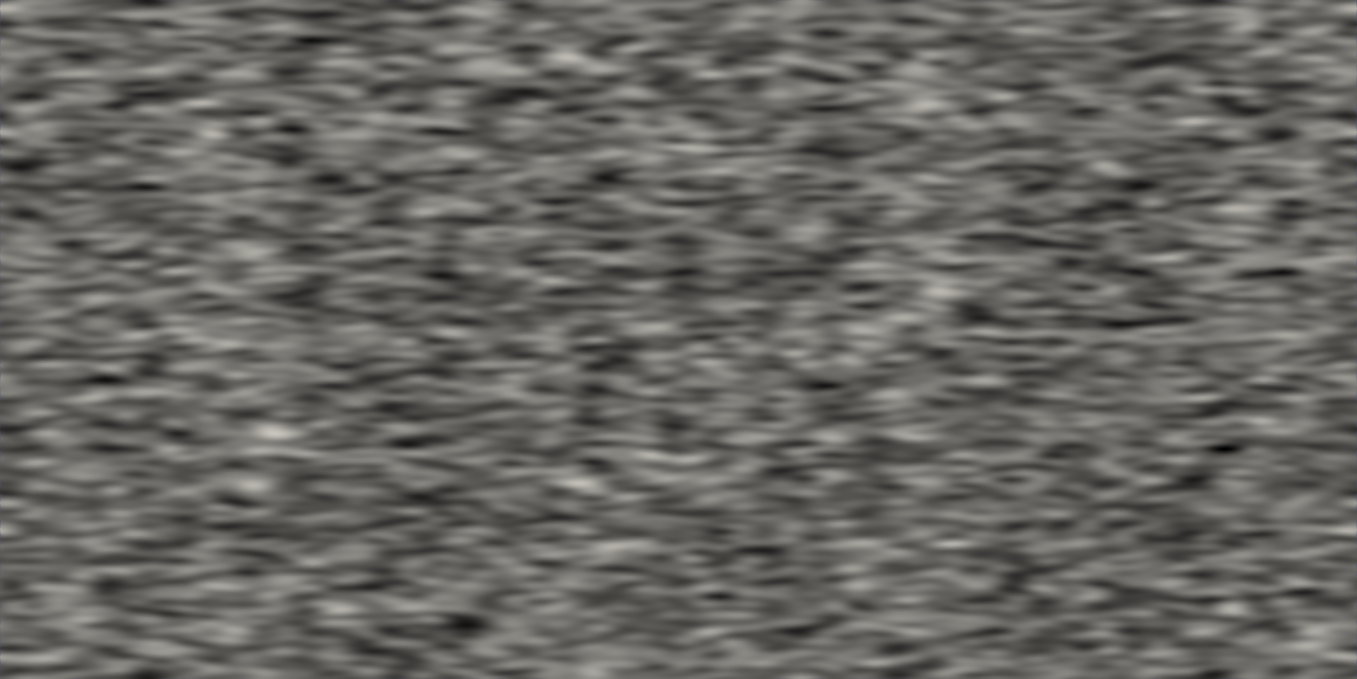}
          \includegraphics[width=0.49\linewidth]{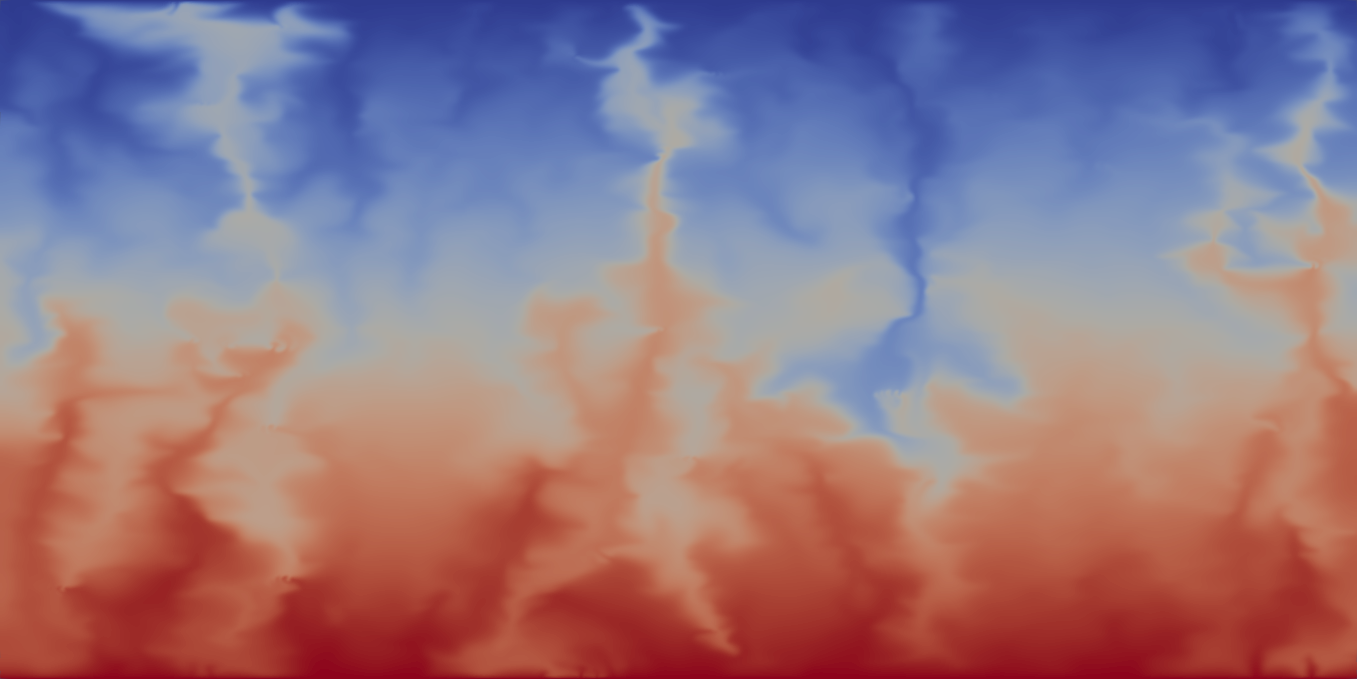}
          \caption{Log-permeability field (left) and temperature field at dimensionless time 5.}\label{fig:HRL}
        \end{figure}

\begin{figure}[!hbtp]
          \includegraphics[width=0.49\linewidth]{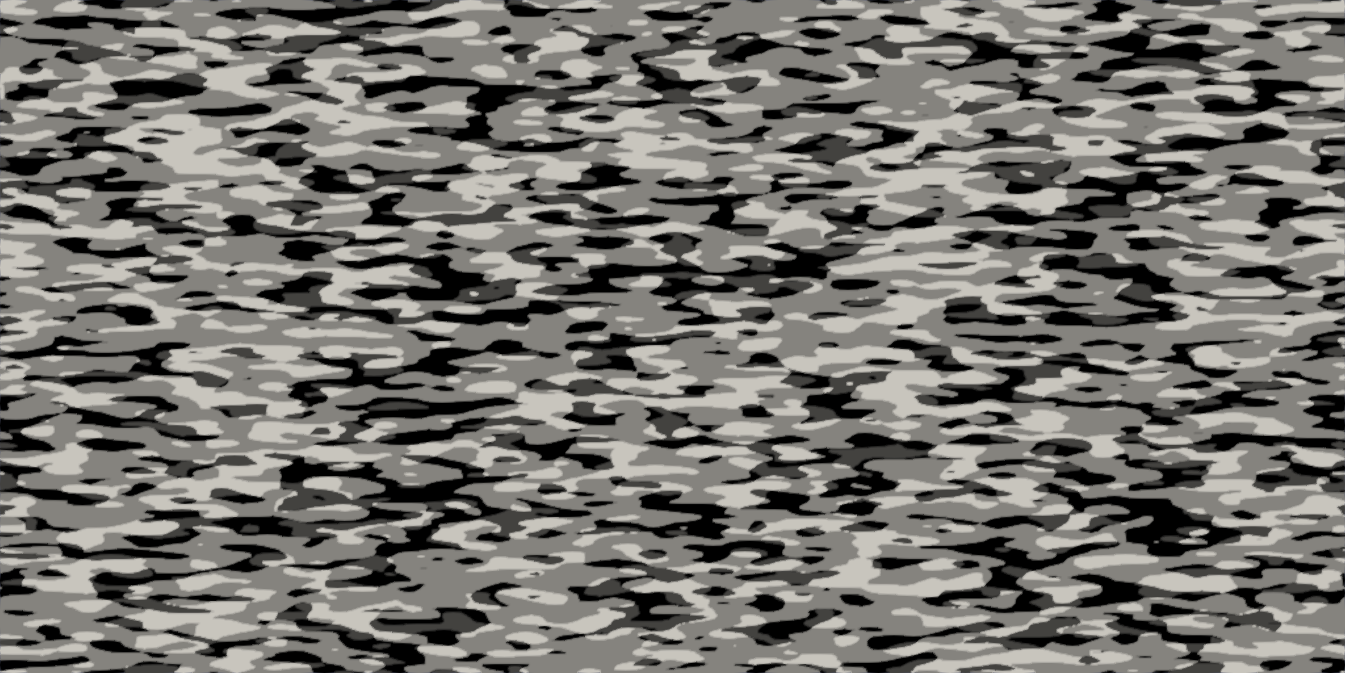}
          \includegraphics[width=0.49\linewidth]{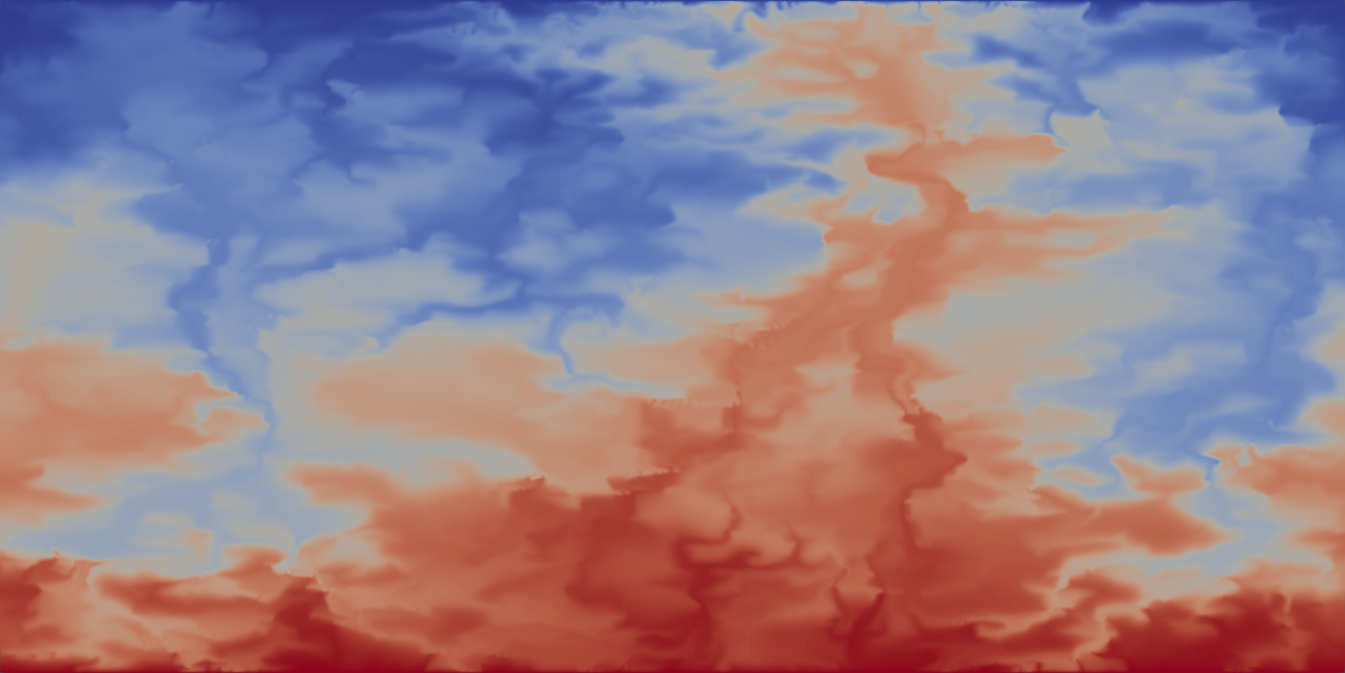}
          \caption{Truncated log-permeability field (left) and temperature solution at dimensionless time 5}\label{fig:HRL-Trun}
        \end{figure}

        For illustration purposes we consider a $2 \times 1$ domain discretised using $1024 \times 512$ cells. Temperature 1 and 0 is prescribed at the bottom and top boundaries respectively and zero temperature gradient is prescribed at the lateral boundaries. Initially temperature varies smoothly with depth. The initial time step is set to 0.1 and the maximum Courant number to 0.5. Two heterogenous cases are solved. First, a log-normally distributed permeability field with $\sigma^{2} = 2$ and correlation lengths  $l_{x} = 0.2; l_{z} = 0.05$ (Figure~\ref{fig:HRL} left). Second, a truncated permeability field with thresholds $2.7\cdot 10^{-1}, 58.9, 1.83 \cdot 10^{-2}$, and 3.99) chosen from the previous permeability distribution (Figure~\ref{fig:HRL-Trun} left). Flow and transport parameters are chosen so that $Ra = 10^{4}$. Solutions for temperature at time 5 are shown in the right panels of figures~\ref{fig:HRL} and \ref{fig:HRL-Trun}.
        
It is important to notice, compared to the classical HRL problem with constant permeability, that the heterogeneities tend to stabilise the flow. In fact, both results in \cref{fig:HRL,fig:HRL-Trun} are steady state, although for other statistically equivalent realisations it has been observed that small local fluctuations could still happen for variance of the log-permeability up to 2. The discontinuities in the permeability field creates, as expected, more defined structures and further stabilise the flow. 
	\subsection{Dual-porosity and discontinuous permeabilities}\label{subsec:dualpor}
        
\begin{figure}[!hbtp]
          \includegraphics[width=0.49\linewidth]{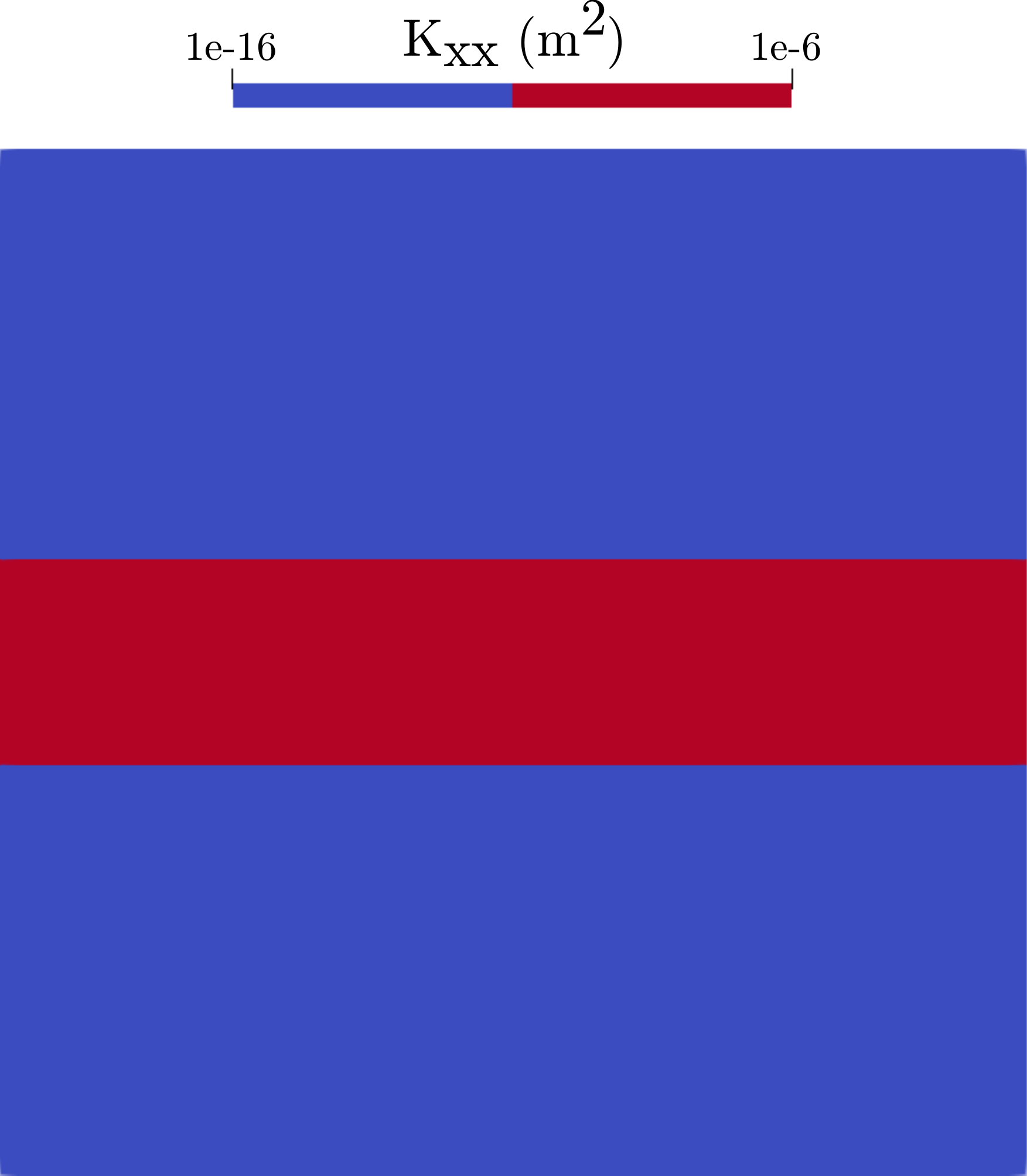}
          \includegraphics[width=0.49\linewidth]{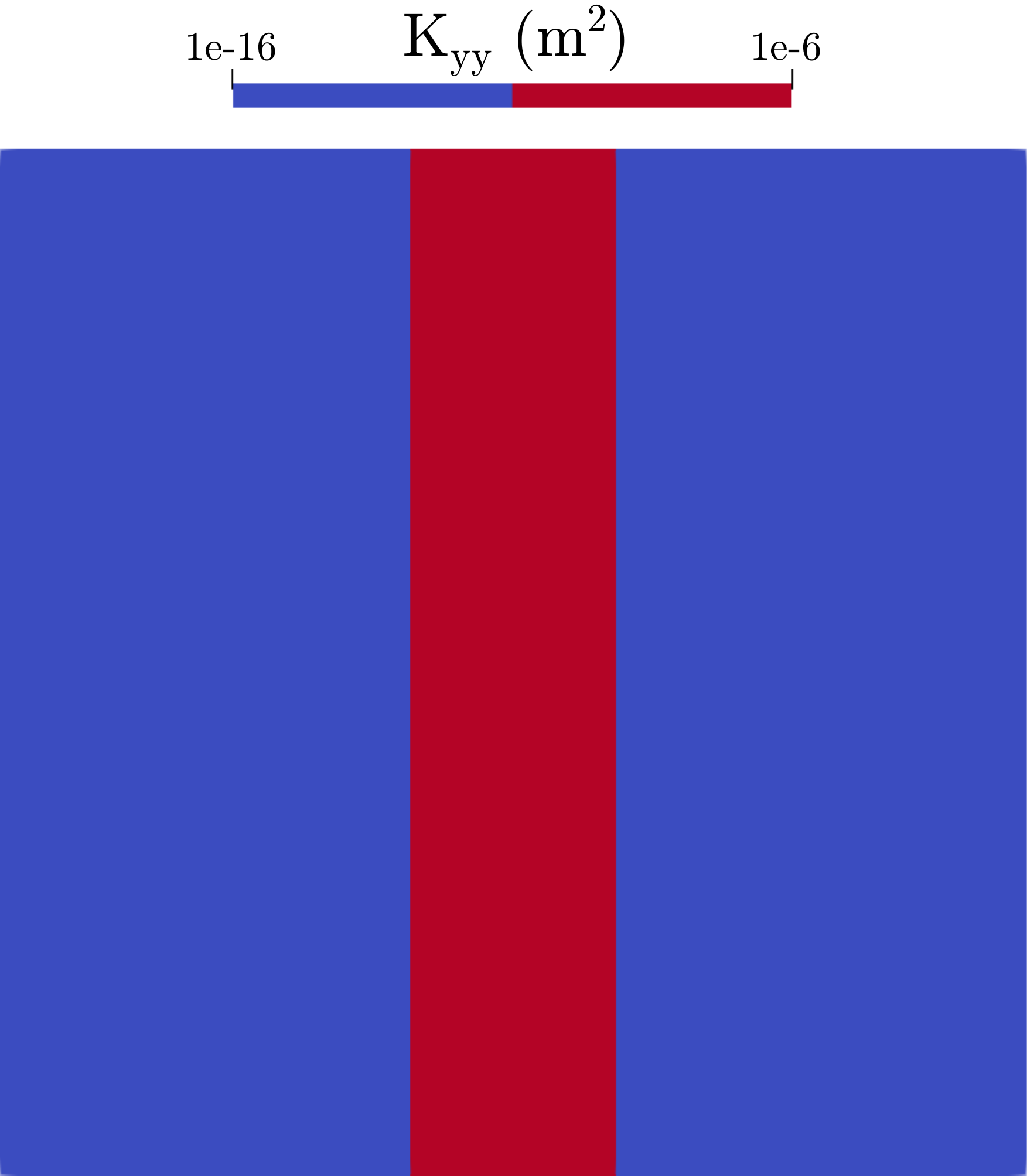}
          \caption{Five-spots injection problem. The low permeability matrix (blue) is  traversed by two anisotropic high permeable fractures (red).}\label{fig:fiveSpot1}
        \end{figure}

        Dual porosity models are a convenient representation of fractures porous media \citep{Barenblatt1960}. In the following example we consider a quarter of a five spot geometry \citep{Simmons1959} in which a low permeability matrix is traversed by two fractures with high anisotropic permeability. Although this model is suitable to describe complex relatively homogeneous networks of fractures, for demonstration purposes, we focus here on a testcase where the fractures are localised in a cross-like structure. Therefore, the fractures permeability tensors have only one non-zero component corresponding to the orientation of the fracture. That is, only the $K_{yy}$ component is non-zero for the vertical fracture and the $K_{xx}$ component for the horizontal one. An injection takes place in the lower left corner and an extraction of the same magnitude on the upper right corner (Figure~\ref{fig:fiveSpot1}). The underlying matrix permeability is isotropic and constant and set to an intermediate value $K_{mat}=10^{-11}$ m$^2$. The matrix porosity is chosen to be $0.3$ for the matrix, dropping to $10^{-5}$ in the regions where fractures are present, while the fracture porosity is close to one ($0.99$) in the fracture-dominated regions and close to zero ($10^{-5}$ s$^{-1}$) elsewhere. The linear transfer coefficient $\tau$ is constant and equal to $10^{-5}$. A Cartesian mesh of $100\times100$ is used.
        
        \begin{figure}[!hbtp]
          \centering
          \includegraphics[width=0.49\linewidth]{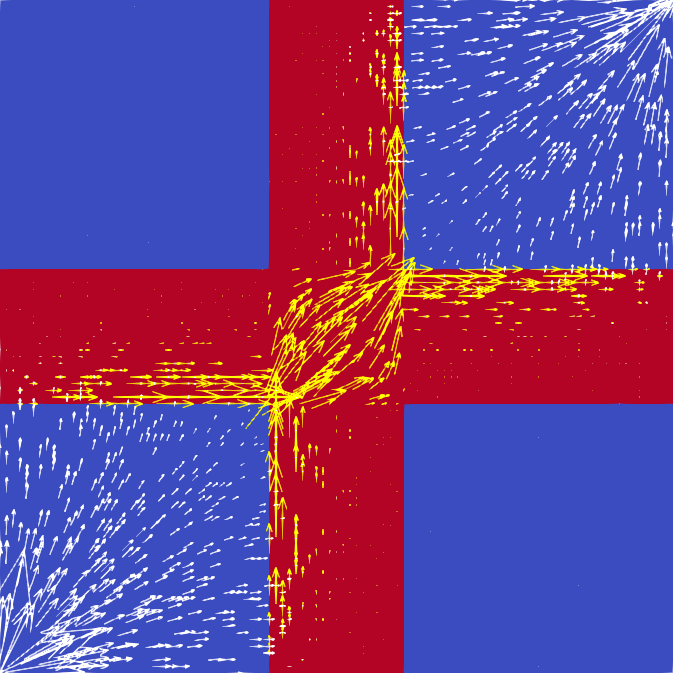}
          \caption{Velocity field for the quater of a five-spot injection problem. Injected water moves preferentially through the matrix (white arrows) until it reaches the high permeability factures (yellow arrows) which focus the flow through the domain. Water leaves the system through the extraction point in the upper right corner.}\label{fig:fiveSpot-vel}
        \end{figure}

The resulting system has a very strong coupling as the flow can only traverse the domain by communicating through the fracture system. This means that the iterative solution of the coupled system via the sequential solutions of matrix and fracture pressure would be very slow (up to thousands of iterations). Thanks to the Schur-based decoupling scheme implemented (the interested reader is referred to \cite{splitting}), the convergence is achieved in only four iterations.
\Cref{fig:fiveSpot-vel} show the resulting velocity field with the sharp transitions between matrix-dominated and fracture-dominated region. For this type of system, where fractures are very localised, the dual-porosity model assumptions breaks down, and similar results could have been obtained with a single highly heterogeneous permeability field. However, the aim here is to test the robustness of the method for solving highly heterogeneous dual-porosity permeability fields. For smoother transitions and larger-scale testcases, where one has a non-negligible fracture and matrix porosity (and permeability) everywhere, the dual-porosity model assumptions would be instead verified. Our geostatistical generation of the fields could be used to produce more realistic scenarios but with the additional difficulty of prescribing a relation or correlation between, not only porosity and permeability of the matrix (as commonly done), but also transfer coefficient and fracture porosity and permeability.

	\section{Conclusions}
In this work we have presented the general-purpose open-source framework \texttt{SECUReFoam} for the computational modelling of various porous media flow and transport problems. An important element of this work is the combination of geostatistical tools with partial differential equations solvers. We focused here on a very general yet simple approach to generate discontinuous random fields (for porosity, permeability and other properties), namely the the truncated pluri-Gaussian simulation. We have limited our attention to unconditioned fields although it is possible to extend this to assimilate real-data by conditioning the random fields. The numerical implementation of continuous and truncated random fields is  based on the explicit spectral representation evaluated in each mesh points rather than Discrete Fourier Transform, therefore independent from the space discretisation and making it suitable for the use with non-structured, non-orthogonal, and locally refined meshes, which are very common in geosciences. Future work will include the conditioning on real data (kriging) and the extension to arbitrary anisotropic correlation matrices including rotations, or more general transformation matrices, allowing therefore more flexibility in the orientation of the facies.

Several mathematical models are implemented, focusing on Darcy's flow, including dual-porosity media, and variable-density flow with different correlation models available for viscosity, density and other parameters. The computational framework also include more advanced physics, such as poro-mechanics, unsaturated flows, and phase-field methods but, due to the physical complexity of these models, these will be the focus of later works. All equations are solved sequentially within a fully implicit Picard-type iteration which deals with non-linearities and coupling between the equations. Operator-based preconditioning and relaxation is used to improve the coupling between the equations. Future work will include the implementation of Newton and Krylov-based iterations for coupled non-linear operators.

We test the framework with well-known benchmark problems, namely the Henry problem for variable-density Darcy flow, the Horton-Rodgers-Lapwood testcase for Rayleigh-Benard instability, a three-dimensional pipe flow for viscous fingering, and a quarter five-spot problem for the dual-porosity model. We demonstrate how the computational framework is robust to highly heterogeneous media, instabilities and coupled problems. All solvers and testcases are available open-source online \cite{federico_municchi_2022_6958098,tutorials}.

\section*{Declarations}
The authors declare that they have no conflicts of interest.

This work has been funded by the following grants:
\begin{itemize}
\item MI, EP, and FM have been supported by the European Union’s Horizon 2020 research and innovation programme, grant agreement number 764531, ”SECURe – Subsurface Evaluation of Carbon capture and storage and Unconventional risks”.
        
\item JJH acknowledges support from Grants CEX2018-000794-S and PID2019-106887GB-C31 funded by AEI/MCINN (10.13039/501100011033) and grant  RYC-2017-22300 funded by MCIN/AEI (10.13039/501100011033) and the European Social Fund.
\end{itemize}

\begin{appendices}

\section{Covariance functions and spectra}
\label{app:covariances}

Covariance functions compute the covariance value $\gamma({r})$ between a pair variables located at points separated by the distance ${r}$.
We focus here on the following cases:

\paragraph{Gaussian}
\begin{equation} \label{eq:gauSpe}
	\gamma({r}; \lambda) = \sigma^2 \left( 1 - e^{-\frac{\pi}{4} \left(\frac{{r}}{\lambda}\right)^2} \right)
	\quad\text{and}\quad
	S(\mathbf{k}; \lambda) = \sigma^2 \left( \frac{\lambda}{\pi} \right)^d e^{- \frac{1}{\pi}(\mathbf{k} \lambda)^2}
\end{equation}

\paragraph{Exponential}
\begin{equation} \label{eq:expSpe}
	\gamma({r}; \lambda) = \sigma^2 \left( 1 - e^{-\frac{|{r}|}{\lambda}} \right)
	\quad\text{and}\quad
	S(\mathbf{k}; \lambda) = \sigma^2 \lambda^d \frac{\Gamma \left( \frac{d+1}{2} \right)}{\left( \pi \left(1 + \mathbf{k}^2 \lambda^2 \right) \right)^{\frac{d+1}{2}}}
\end{equation}

\paragraph{Matérn} \cite{wackernagel1998examples}
\begin{equation} \label{eq:matSpe}
\begin{split}
	\gamma({r}; \nu, \lambda) & =
	\sigma^2 \left( 1 -  \frac{2^{1-\nu}}{\Gamma(\nu)} \left( \sqrt{2 \nu} \frac{{r}}{\lambda} \right)^{\nu} K_{\nu} \left( \sqrt{2 \nu} \frac{{r}}{\lambda} \right) \right)
	\\
	S(\mathbf{k}; \nu, \lambda) & = \sigma^2 \lambda^{d} \frac{\Gamma(\nu + \frac{d}{2})(2 \nu)^{\nu}}
	{\Gamma(\nu)\pi^{\frac{d+1}{2}}\left( {2 \nu} + \lambda^2 \mathbf{k}^2 \right)^{\nu + \frac{d}{2} }}
\end{split}
\end{equation}
where $K_{\nu}$ is the modified Bessel function of the second kind.

\paragraph{Spherical} \cite{wackernagel1998examples}
\begin{equation}
  \gamma({r}; \lambda) = \sigma^2
  \begin{cases}
	  1 - 1.5 \frac{{r}}{\lambda} + 0.5 \left(\frac{{r}}{\lambda}\right)^3, & \text{ if } {r}<\lambda \\
	  0, & \text{ otherwise}
  \end{cases}
\end{equation}


%

\section{Random field generator options}
\label{app:setRandomFieldDict}
The users can tune their models on two different levels: by setting pre-defined key words and numerical fields in a dictionary or editing the source code for more substantial changes. The tuning options available in the dictionary are:
\begin{itemize}
	\item \texttt{type} of GRF: \texttt{truncated} for ordinate geological formations (i.e. facies are self-embedded) or \texttt{bitruncated} for domains where facies show more complicate pattern. The difference at computational level is that the former generates a single GRF while the latter produces two GRFs;
	\item \texttt{correlation} function: \texttt{gaussian}, \texttt{exponential} or \texttt{matern};
	\item fields and correlation parameters:
	\begin{itemize}
		\item \texttt{Kmean} and \texttt{Ksigma}: mean and variance of the Gaussian fields;
		\item \texttt{Lcorr}: correlation lengths;
		\item \texttt{$\rho$} and \texttt{$\nu$}: just for \texttt{matern};
		\item \texttt{nfreq}: correlation spectrum parameters;
	\end{itemize}
	\item \texttt{disableY} and/or \texttt{disableZ}: they provide control over the number of dimensions of the field. When they are \texttt{false} the generated TGS is 3D;
	\item \texttt{periodic}: it controls the periodicity structure of the domain;
	\item \texttt{values}: float numbers taken by the facies;
	\item \texttt{thresholds}: float array which sets the boundaries of the of the facies on the chosen Gaussian distribution;
	\item \texttt{thresholds2}: same as \texttt{thresholds} but for the second Gaussian distribution, in case \texttt{bitruncated} option is selected.
\end{itemize}

\section{Density and viscosity models}
\label{app:density_viscosity}

Pre-defined options in the library that can be selected independently for $\rho(c)$ and $\mu(c)$. The models are:
\begin{description}
  \item[\texttt{constant}:] A costant value is used for the fluid property.
  \item[\texttt{linear}:] Fluid property changes linearly with concentration as $f(c) = f_{0} + f_{1}c$.
  \item[\texttt{exponential}:] An exponential relation of the form $f(c) = f_{0}e^{Rc}$ is used for density or viscosity.
  \item[\texttt{tabulated}:] Fluid properties are interpolated from tabulated data $(c, f)$ allowing the use of more complex functions or experimental data. The interpolation is done using the \texttt{nonUniformTable} class of \texttt{OpenFoam} \texttt{thermophysicalFunctions}.
\end{description}
\end{appendices}

	\bibliography{references}

\end{document}